\documentclass{bryant-article}

\usepackage{color}
\usepackage{enumerate}

\newcommand{\bu}{\mathbf{u}}
\newcommand{\ueq}{\mathrel{\dot=}}
\newcommand{\lk}{\operatorname{lk}}
\newcommand{\ab}{\operatorname{ab}}
\newcommand{\adj}{\operatorname{adj}}
\newcommand{\spn}{\operatorname{span}}
\newcommand{\Fox}[2]{\frac{\partial #1}{\partial #2}}
\DeclareMathOperator{\ord}{mod}

\title[Alexander polynomials of twisted torus links]{Multivariable Alexander polynomials of twisted torus links via Dehn surgery}

\author{Birch Bryant}
\address{University of North Alabama}
\email{bbryant3@una.edu}

\subjclass[2020]{57K10, 57K14}
\keywords{Alexander polynomial, twisted torus knot, Dehn surgery, Fox calculus, Burau representation}

\begin{document}
	
	\begin{abstract}
		Let $T(p,q;r,s)$ be the twisted torus link obtained from the torus link $T(p,q)$ by adding $s$ full twists to $r$ adjacent strands, and let $C$ be an unknotted circle encircling those strands. For all $p$, $q$, and $r$ we compute the multivariable Alexander polynomial of $T(p,q)\cup C$ in closed form. A relation of Silver and Williams between the Alexander polynomials of a link and of its $1/n$ Dehn fillings then gives a closed formula for the multivariable Alexander polynomial of every twisted torus link $T(p,q;r,s)$, answering a question of Adnan and Park. When $\gcd(p,q)=1$, so that $T(p,q;r,s)$ is a knot, the formula is a sum of $pr$ monomials with coefficient $+1$, divided by cyclotomic factors, and its exponents are governed by an elementary counting function on the orbit of $q$ modulo $p$. Because the numerator has no cancellation, we are able to give the degree for all $s$ and simple non-fiberedness and non-L-space criteria for infinite families of twisted torus knots.
	\end{abstract}
	
	\maketitle

	\section{Introduction}\label{sec:intro}

	Fix integers $1<q<p$, $1\le r\le p-1$, and $s\in\mathbb Z$. Let $\delta=\sigma_1\sigma_2\cdots\sigma_{p-1}\in B_p$ be a word in the braid group on $p$ strands, so that the closure of $\delta^q$ is the torus link $T(p,q)$, and let $\theta_r=(\sigma_1\cdots\sigma_{r-1})^r$ be the full twist on the first $r$ strands. The \emph{twisted torus link} is the braid closure
	\[
	T(p,q;r,s) \;=\; \widehat{\;\delta^{q}\,(\theta_r)^{s}\;}.
	\]
	Equivalently, $T(p,q;r,s)$ is obtained from $T(p,q)$ by $-1/s$ surgery on an unknotted circle $C$ that bounds a disc meeting $T(p,q)$ in $r$ adjacent strands. It has $g=\gcd(p,q)$ components; when $g=1$ it is the \emph{twisted torus knot} introduced by Dean~\cite{Dean} in the study of Seifert-fibered surgeries.
	
	The Alexander polynomials of twisted torus knots have been studied by several authors. Morton~\cite{Morton} computed $\Delta_{T(p,q;2,s)}$ for all $s$, using the multivariable Alexander polynomial of a three-component link formed by the knot and the cores of the two solid tori of the genus-one Heegaard splitting. Adnan and Park~\cite{AdnanPark} recently gave a formula for general $r$ using Fox calculus on a presentation coming from a genus-two Heegaard surface. They determined the degree for $s>0$ and for $s<0$ with $r|s|<q$. They noted that the general negative case is more difficult, and asked for the Alexander polynomials of twisted torus links~\cite[Question~6]{AdnanPark}. 
	
	In this paper we take a Dehn surgery approach. The primary tool is a relation given in Silver and Williams~\cite[Remark~3.3]{SilverWilliams}. If $L=L_1\cup\dots\cup L_d$ is a link, $\lambda_i=\lk(L_i,L_d)$, and $L(n)$ denotes the result of $1/n$ surgery on $L_d$, then
	\begin{equation}\label{eq:SW}
		\Delta_L\Bigl(u_1,\dots,u_{d-1},\prod_{i<d}u_i^{-n\lambda_i}\Bigr)
		\ueq
		\begin{cases}
			\bigl(\prod_{i<d}u_i^{\lambda_i}-1\bigr)\,\Delta_{L(n)}(u_1,\dots,u_{d-1}), & d\ge3,\\[3pt]
			\dfrac{u_1^{\lambda_1}-1}{u_1-1}\,\Delta_{L(n)}(u_1), & d=2.
		\end{cases}
	\end{equation}
	Here $\ueq$ denotes equality up to multiplication by a unit $\pm u_1^{e_1}\cdots u_{d-1}^{e_{d-1}}$. Thus a closed form for the multivariable Alexander polynomial of $T(p,q)\cup C$ gives the Alexander polynomials of the whole family $\{T(p,q;r,s)\}_{s\in\mathbb Z}$ at once. 
	
	This relation strikes the author as an underutilized tool for Alexander polynomial computation for families of knots; the limiting factor being derivation of $\Delta_L$. It is to this task we set ourselves for the links $T(p,q)\cup C$.
	
	\subsection*{Notation}
	Let $g=\gcd(p,q)$, $p'=p/g$, and $q'=q/g$. The components of $T(p,q)$ are $L_0,\dots,L_{g-1}$, where strand $j$ of $\delta^q$ lies on $L_{c(j)}$ with $c(j)=(j-1)\bmod g$ (Section~\ref{sec:colors}). Let $u_0,\dots,u_{g-1}$ be the variables of the components, write $\bu=(u_0,\dots,u_{g-1})$, and let $v$ be the variable of $C$. Set
	\[
	U=u_0u_1\cdots u_{g-1},\qquad V=u_0^{\lambda_0}\cdots u_{g-1}^{\lambda_{g-1}},\qquad G_i=u_{c(1)}u_{c(2)}\cdots u_{c(i)},
	\]
	where
	\[
	\lambda_c=\#\{1\le j\le r: c(j)=c\}.
	\]
	Here $\lambda_c=\lk(K_c,C)$ and $G_0=1$. For $0\le i\le r-1$ and $k\ge0$ define the \emph{counting function}
	\begin{equation}\label{eq:N}
		N_i(k)\;=\;\#\bigl\{\,m\in\{1,\dots,k\}\;:\;(mq+i)\ \ord\ p\;<\;r\,\bigr\}.
	\end{equation}
	That is, $N_i(k)$ counts how many of the first $k$ points of the orbit $q,2q,3q,\dots\pmod p$ land in the cyclic window $\{-i,\dots,r-1-i\}$ of length $r$. Finally, for $0\le c\le g-1$ let
	\[
	S_c(\bu,v)=\sum_{n=0}^{p'-1}\ \sum_{\substack{0\le i\le r-1\\ i\equiv c\ (\mathrm{mod}\ g)}} U^{nq'}\,G_i\,v^{N_i(n)}.
	\]
	When $g=1$ there is one variable $u$, and $U=u$, $V=u^r$, $G_i=u^i$, and $S_0(u,v)=\sum_{k=0}^{p-1}\sum_{i=0}^{r-1}u^{kq+i}v^{N_i(k)}$.
	
	\subsection*{Main results} We find for the augmented torus knot $T(p,q)\cup C$ a formula for the $g+1$ variable Alexander polynomial.
	
	\begin{theorem}\label{thm:A}
		Orient $C$ so that $\lk(K_c,C)=\lambda_c\ge0$ for all $c$. Then
		\[
		\Delta_{T(p,q)\cup C}(\bu,v)\;\ueq\;\frac{\displaystyle\sum_{c=0}^{g-1}(1-u_c)\prod_{c'\neq c}\bigl(1-U^{p'q'}v^{\lambda_{c'}}\bigr)\,S_c(\bu,v)}{U^{p'}-1}.
		\]
		In particular, if $\gcd(p,q)=1$, then
		\[
		\Delta_{T(p,q)\cup C}(u,v)\;\ueq\;\frac{1-u}{1-u^{p}}\;\sum_{k=0}^{p-1}\sum_{i=0}^{r-1} u^{\,kq+i}\,v^{\,N_i(k)}.
		\]
	\end{theorem}
	
	Filling $C$ by~\eqref{eq:SW} gives the twisted torus links and knots. For links ($g\ge2$) the link $T(p,q)\cup C$ has $g+1\ge3$ components, and we use the case $d\ge3$.
	
	\begin{theorem}\label{thm:links}
		Let $g=\gcd(p,q)\ge2$. For all $s\in\mathbb Z$,
		\[
		\Delta_{T(p,q;r,s)}(\bu)\;\ueq\;
		\frac{\displaystyle\sum_{c=0}^{g-1}(1-u_c)\prod_{c'\neq c}\bigl(1-U^{p'q'}V^{s\lambda_{c'}}\bigr)\,S_c\bigl(\bu,V^{s}\bigr)}{(U^{p'}-1)(V-1)}.
		\]
	\end{theorem}
	
	For knots ($g=1$) we use the case $d=2$, and the formula takes a particularly simple form.
	
	\begin{theorem}\label{thm:knots}
		Let $\gcd(p,q)=1$. For all $s\in\mathbb Z$,
		\[
		\Delta_{T(p,q;r,s)}(t)\;\ueq\;\frac{(1-t)^{2}}{(1-t^{p})(1-t^{r})}\sum_{k=0}^{p-1}\sum_{i=0}^{r-1} t^{kq+i+rs\,N_i(k)}.
		\]

	\end{theorem}
	
	Every coefficient of the double sum is positive, so the extreme degrees of $\Delta_{T(p,q;r,s)}$ cannot cancel. Together with a symmetry of $N_i(k)$ (Lemma~\ref{lem:pairing}), this gives the degree for every $s$. We write $\spn f$ for the difference between the largest and smallest exponents of a Laurent polynomial $f$, so that $\deg \Delta_K=\tfrac12\spn\Delta_K$ in the symmetric normalization. Let $E(k,i)= kq+i+rsN_i(k)$.
	
	\begin{theorem}\label{thm:degree}
		Let $\gcd(p,q)=1$. For all $s\in\mathbb Z$,
		\[
		\spn \Delta_{T(p,q;r,s)} \;=\; 2\max_{k,i} E(k,i)\;-\;(p-1)q\;-\;(r-1)(1+rs)\;-\;(p+r-2).
		\]
	\end{theorem}
	
	\begin{corollary}\label{cor:fiber}
		Let $\gcd(p,q)=1$. The leading coefficient of $\Delta_{T(p,q;r,s)}$ equals, up to sign, the number of pairs $(k,i)$ at which $E(k,i)$ is maximal. In particular, if this maximum is attained more than once, then $T(p,q;r,s)$ is not fibered.
	\end{corollary}
	
	For $s>0$ the knot $T(p,q;r,s)$ is the closure of a positive braid, so it is fibered, and Corollary~\ref{cor:fiber} says nothing new. The interesting case is $s<0$. For example, $\Delta_{T(4,3;2,-2)}=2t^2-3t+2$, so $T(4,3;2,-2)$ is not fibered. Adnan and Park constructed infinite families of non-fibered twisted torus knots from the leading coefficients of their formula~\cite{AdnanPark2}.
	
	By Ozsv\'ath and Szab\'o \cite{OzsvathSzabo}, the nonzero coefficients of the Alexander polynomial of an L-space knot are $\pm 1$ and alternate sign. For $p \equiv \pm 1 \mod q$ and $s>0$, Vafaee \cite{Vafaee} determined exactly which $T(p,q;r,s)$ are L-space knots. Reading the Alexander polynomial from Theorem \ref{thm:knots} gives obstructions with $s<0$.
	
	\begin{theorem}\label{thm:lspace}
		Let $\gcd(p,q)=1$ and let $p=nq-1$ with $n\geq 2$, let $2\leq r \leq q-1$, and let $s<0$ with $r|s| \geq p+3$. Then, up to a unit, $$\Delta_{T(p,q;r,s)}(t)=t^D - 2t^{D-1}+\cdots.$$In particular, $\Delta_{T(p,q;r,s)}$ is monic, but $T(p,q;r,s)$ is not an L-space knot.
	\end{theorem}
	
	\subsection*{Outline}
	In Section~\ref{sec:prelim} background is given on Fox calculus, the Alexander polynomial, and braids. Section~\ref{sec:group} computes the group of the augmented link with braid axis $T(p,q)\cup A \cup C$. In Section~\ref{sec:aug}, the formula for the Alexander polynomial of $T(p,q)\cup A\cup C$ is proven and Theorem~\ref{thm:A} is then deduced. Section~\ref{sec:links} gives the proofs of Theorem~\ref{thm:links} for twisted torus links and the companion result Theorem~\ref{thm:knots} for knots. In Section 6, we explore the properties of the counting function $N_i(k)$ and prove Theorem \ref{thm:degree}, Corollary \ref{cor:fiber}, and Theorem \ref{thm:lspace}. We close by giving further questions in Section 7.

	\section{Background}\label{sec:prelim}

	\subsection{Fox calculus}\label{sec:fox}
	Let $F=F(x_1,\dots,x_n)$ be free. The Fox derivatives \\$\partial/\partial x_j\colon F\to\mathbb Z F$ are characterized by $\Fox{x_i}{x_j}=\delta_{ij}$ and $\Fox{(uv)}{x_j}=\Fox{u}{x_j}+u\Fox{v}{x_j}$. They satisfy the \emph{fundamental formula}~\cite{Fox53}
	\begin{equation}\label{eq:fundamental}
		\sum_{j=1}^n \Fox{w}{x_j}\,(x_j-1)=w-1\qquad(w\in F).
	\end{equation}
	Let $\Lambda$ be a Laurent polynomial ring $\mathbb Z[t_1^{\pm1},\dots,t_d^{\pm1}]$, and let $\ab\colon\mathbb Z F\to\Lambda$ be the ring homomorphism sending each $x_j$ to one of the variables $t_c$ (in our applications, the meridian class of the component that $x_j$ goes around). For an endomorphism $\phi$ of $F$, the \emph{Fox Jacobian} is the matrix $\bigl(\Fox{\phi(x_i)}{x_j}\bigr)_{i,j}$ over $\mathbb Z F$, and the \emph{abelianized Jacobian} is its image
	\[
	J(\phi)=\Bigl(\ab\,\Fox{\phi(x_i)}{x_j}\Bigr)_{i,j}\in M_n(\Lambda).
	\]
	
The Fox derivatives satisfy the chain rule $\Fox{(\phi\circ\psi)(x_i)}{x_j}=\sum_k\phi\bigl(\Fox{\psi(x_i)}{x_k}\bigr)\Fox{\phi(x_k)}{x_j}$. If $\ab\circ\phi=\ab$, applying $\ab$ gives
\[
J(\phi\circ\psi)=J(\psi)\,J(\phi).
\]
The hypothesis $\ab\circ\phi=\ab$ holds for every braid automorphism when all $x_i$ are sent to the same variable, but not in general: for links, a braid may move a strand to a different component (see Section~\ref{sec:links}).
	
	\subsection{Alexander polynomials and deficiency-one presentations}
	Let $L$ be an oriented link with $d$ components, and let $\Lambda=\mathbb Z[t_1^{\pm1},\dots,t_d^{\pm1}]$, with $t_c$ the class of a positive meridian of the $c$-th component. For a finite presentation $\langle x_1,\dots,x_n\mid R_1,\dots,R_m\rangle$ of $\pi_1(S^3\setminus L)$, the \emph{Alexander matrix} is $\bigl(\ab\,\Fox{R_i}{x_j}\bigr)$. Its first elementary ideal $E_1$, generated by the $(n-1)$-minors, does not depend on the presentation~\cite[Ch.~VII]{CrowellFox}. For $d\ge 2$ the Alexander polynomial $\Delta_L$ is the generator of the smallest principal ideal containing $E_1$, well defined up to units of $\Lambda$; for knots this convention is modified by a factor $t-1$ \cite{Torres}.
	
	\begin{lemma}[Torres~\cite{Torres}; see also~\cite{CrowellFox,Kawauchi}]\label{lem:FoxTorres}
		Suppose $d\ge2$, and let $\langle x_1,\dots,x_n\mid R_1,\dots,R_{n-1}\rangle$ be a presentation of $\pi_1(S^3\setminus L)$ of deficiency one in which every generator $x_j$ is a meridian, so that $z_j:=\ab(x_j)-1\in\{t_1-1,\dots,t_d-1\}$. Let $A$ be the Alexander matrix and let $A_j$ be $A$ with its $j$-th column deleted. Then $\det A_j\ueq z_j\,\Delta_L$ for every $j$.
	\end{lemma}

	The Alexander matrices in this paper all come from relations of the form $uv^{-1}$ with $\ab(u)=\ab(v)$, and their entries are differences of abelianized Jacobians.
	
	\begin{lemma}\label{lem:alexjac}
		Let $R=uv^{-1}$ with $u,v\in F$ and $\ab(u)=\ab(v)$. Then $\ab\Fox{R}{x_j}=\ab\Fox{u}{x_j}-\ab\Fox{v}{x_j}$. In particular:
		\begin{enumerate}
			\item[(a)] For a braid $\beta\in B_p$ and the closed-braid relations $r_i=x_i\,\beta(x_i)^{-1}$, the Alexander matrix is $I-J(\beta)$.
			\item[(b)] Let $\phi\colon F(x_1,\dots,x_p)\to F(x_1,\dots,x_p,y_1,\dots,y_m)$ be a homomorphism with $\ab\circ\phi=\ab$ on $F(x_1,\dots,x_p)$, and let $R_i=\tau x_i\tau^{-1}\phi(x_i)^{-1}$ in the free group on $x_1,\dots,x_p,y_1,\dots,y_m,\tau$, with $\ab(\tau)=a$. Then the $x$-columns of the rows $R_1,\dots,R_p$ form the matrix $aI-J(\phi)$, where $J(\phi)=\bigl(\ab\Fox{\phi(x_i)}{x_j}\bigr)$, and the $\tau$-column has entries $\ab\Fox{R_i}{\tau}=1-\ab(x_i)$.
		\end{enumerate}
	\end{lemma}
	
	\begin{proof}
		By the product rule, $\Fox{(uv^{-1})}{x_j}=\Fox{u}{x_j}+u\Fox{(v^{-1})}{x_j}=\Fox{u}{x_j}-uv^{-1}\Fox{v}{x_j}$, and $\ab(uv^{-1})=1$. For (a), $\Fox{x_i}{x_j}=\delta_{ij}$. For (b), $\Fox{(\tau x_i\tau^{-1})}{x_j}=\tau\delta_{ij}$, which abelianizes to $a\delta_{ij}$; and $\Fox{(\tau x_i\tau^{-1})}{\tau}=1-\tau x_i\tau^{-1}$, while $\phi(x_i)$ does not involve $\tau$.
	\end{proof}
	
	We also use the Torres condition~\cite{Torres}: if $L=L'\cup\ell$ has $d\ge 3$ components, then
	\begin{equation}\label{eq:torres}
		\Delta_L(t_1,\dots,t_{d-1},1)\ueq\bigl(t_1^{\lk(L_1,\ell)}\cdots t_{d-1}^{\lk(L_{d-1},\ell)}-1\bigr)\,\Delta_{L'}(t_1,\dots,t_{d-1}).
	\end{equation}
	
	\subsection{Braids}\label{sec:braids}
	For a braid $\beta\in B_p$, the \emph{closure} $\hat\beta\subset S^3$ is the link obtained by joining the top endpoint of each strand of $\beta$ to the bottom endpoint in the same position by arcs running around the braid axis. For example, $\widehat{\delta^q}=T(p,q)$.
	
	Let $D_p$ be a disc with $p$ punctures and a base point on its boundary, and let $x_1,\dots,x_p$ be the standard free generators of $\pi_1(D_p)$, with $x_i$ a loop around the $i$-th puncture. A braid $\beta$ determines a homeomorphism of $D_p$ fixing $\partial D_p$, and hence an automorphism of $F(x_1,\dots,x_p)$. We denote this automorphism by the same letter, so $\beta(x_i)$ is the image of $x_i$. On the generators of $B_p$ it is the Artin action
	\[
	\sigma_j\colon\; x_j\mapsto x_jx_{j+1}x_j^{-1},\quad x_{j+1}\mapsto x_j,\quad x_i\mapsto x_i\ (i\ne j,j+1),
	\]
	and products of braids act by composition, $(\beta_1\beta_2)(x_i)=\beta_1\bigl(\beta_2(x_i)\bigr)$. Thus for $\beta=\sigma_{j_1}^{\varepsilon_1}\cdots\sigma_{j_m}^{\varepsilon_m}$ the substitutions are applied from the last letter to the first. With this convention,
	\[
	\pi_1(S^3\setminus\hat\beta)=\langle x_1,\dots,x_p\mid x_i=\beta(x_i)\ (1\le i\le p)\rangle,
	\]
	one relation being redundant~\cite{KasselTuraev}.
	
	\subsubsection*{The colored Burau matrix}
	Fix a \emph{coloring} $c\colon\{1,\dots,p\}\to\{0,\dots,g-1\}$ of the strand positions and variables $u_0,\dots,u_{g-1}$, and let $\ab$ send $x_j\mapsto u_{c(j)}$. Every braid automorphism sends each generator to a conjugate of a generator~\cite{Birman,KasselTuraev}: there is a permutation $\bar\beta$ of $\{1,\dots,p\}$, the \emph{permutation of $\beta$}, together with words $W_1,\dots,W_p$ such that
	\[
	\beta(x_i)=W_i\,x_{\bar\beta(i)}\,W_i^{-1}\qquad(1\le i\le p).
	\]
	Hence $\ab\circ\beta=\ab$ exactly when $c(\bar\beta(i))=c(i)$ for all $i$; we then say that $\beta$ \emph{preserves the coloring}. For such $\beta$, the \emph{colored Burau matrix} is the abelianized Jacobian
	\[
	\mathcal B(\beta)=J(\beta)=\Bigl(\ab\,\Fox{\beta(x_i)}{x_j}\Bigr)_{i,j}.
	\]

	\subsubsection*{The braids $\delta$ and $\theta_r$}
	We record the automorphisms of the two braids from the introduction. For $\delta=\sigma_1\cdots\sigma_{p-1}$ and $i<p$, the letters $\sigma_{p-1},\dots,\sigma_{i+1}$ fix $x_i$, then $\sigma_i$ gives $x_ix_{i+1}x_i^{-1}$, and $\sigma_{i-1},\dots,\sigma_1$ lower the index of the outer letter to $1$. For $i=p$, the letters lower $x_p$ to $x_1$. Hence
	\begin{equation}\label{eq:delta}
		\delta(x_i)=x_1x_{i+1}x_1^{-1}\quad(1\le i<p),\qquad \delta(x_p)=x_1 .
	\end{equation}
	Geometrically, $\delta$ rotates the punctures cyclically by one position. For the full twist $\theta_r=(\sigma_1\cdots\sigma_{r-1})^r$ on the first $r$ strands, let $w_r=x_1x_2\cdots x_r$. Then
	\begin{equation}\label{eq:fulltwist}
		\theta_r(x_i)=\begin{cases} w_r\,x_i\,w_r^{-1}, & i\le r,\\ x_i, & i>r.\end{cases}
	\end{equation}
	Indeed, $\theta_r$ rotates a subdisc $E\subset D_p$ containing the first $r$ punctures through $360^\circ$ and is the identity outside $E$. A loop around a puncture in $E$ is therefore conjugated by the boundary loop of $E$, which represents $w_r$.

	\section{The group of the augmented link}\label{sec:group}

		\begin{figure}[ht]
		\centering
		\includegraphics[width=3in]{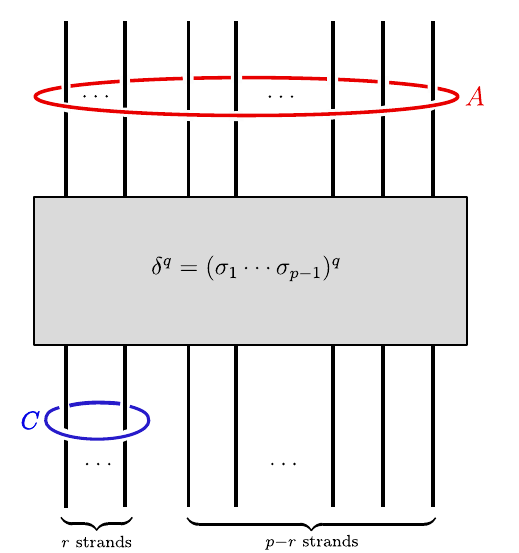}
		\caption{The link $T(p,q)\cup A\cup C$. The circle $C$ encircles the first $r$ strands at a single level of the braid, and $A$ is the braid axis. }
		\label{fig:link}
	\end{figure}
	
	Let $\beta\in B_p$, let $A$ be the axis of the closed braid $\hat\beta$, and let $C\subset S^3\setminus(\hat\beta\cup A)$ bound a disc $E$ that lies in a page of the open book with binding $A$ and meets $\hat\beta$ transversely in the $r$ adjacent points $1,\dots,r$ (Figure~\ref{fig:link}). Define the homomorphism
	\[
	\psi_c\colon F(x_1,\dots,x_p)\to F(x_1,\dots,x_p,c),\qquad
	\psi_c(x_i)=\begin{cases}c\,x_i\,c^{-1},& i\le r,\\ x_i,& i>r,\end{cases}
	\]
	which conjugates the first $r$ generators by $c$.
	
	\begin{lemma}\label{lem:presentation}
		With $w=x_1\cdots x_r$, the group $\pi_1\bigl(S^3\setminus(\hat\beta\cup A\cup C)\bigr)$ has the presentation with generators $x_1,\dots,x_p,c,\tau$ and relations
		\[
		\tau x_i\tau^{-1}=\psi_c\bigl(\beta(x_i)\bigr)\quad(1\le i\le p),\qquad c\,w\,c^{-1}=w,
		\]
		where $x_i$, $c$, and $\tau$ are meridians of $\hat\beta$, $C$, and $A$.
	\end{lemma}
	
	\begin{proof}
		 The complement $V$ of an open tubular neighborhood of $A$ is a solid torus $D\times S^1$ whose meridian discs are the pages of the open book. The closed braid meets each page transversely in $p$ points. Hence $X=V\setminus N(\hat\beta)$ is the mapping torus $D_p\times[0,1]/(y,0)\sim(h(y),1)$ of the $p$-punctured disc, with monodromy $h$ fixing $\partial D$. Choose a base point $*\in\partial D$ and the standard generators $x_1,\dots,x_p$ of $\pi_1(D_p,*)$, so that $h_*=\beta$. The arc $\nu=\{*\}\times[0,1]$ closes up to a loop $\tau$ on $\partial V$ that meets each page once, so $\tau$ is a meridian of $A$.
		
		After an isotopy in the complement of $\hat\beta\cup A$, we may assume $C=\partial E\times\{\tfrac12\}$, where $E\subset D$ is a disc containing exactly the punctures $1,\dots,r$ and $*\notin E$. Because these punctures are adjacent, the generators can be chosen so that $\partial E$, joined to $*$ by an arc in $D\setminus E$, represents $w=x_1\cdots x_r$.
		
		Now, let $Y=(D_p\times[0,1])\setminus N(C)$, and let $Y_\pm$ be its parts of height $\ge\tfrac12$ and $\le\tfrac12$. Each deformation retracts onto a copy of $D_p$. Thus $\pi_1(Y_-)$ is free on $x_1,\dots,x_p$ and $\pi_1(Y_+)$ is free on $x_1',\dots,x_p'$, the latter transported to height $0$ along $\nu$. The intersection $Y_+\cap Y_-$ is the punctured page with $N(C)$ removed, and it has two components:
		\begin{itemize}
			\item $S_{\rm out}=(D_p\setminus E)\times\{\tfrac12\}$, whose fundamental group is free on $x_{r+1},\dots,x_p$ and $w$;
			\item $S_{\rm in}=(E\cap D_p)\times\{\tfrac12\}$, a disc with $r$ holes, based at a point $e$.
		\end{itemize}
		Choose an arc $\alpha\subset D$ from $*$ to $e$, and let $\alpha_\pm$ be its push-offs into $Y_\pm$. Set $c=\alpha_+\alpha_-^{-1}$. This loop passes once through $E$, so it is a meridian of $C$. The van Kampen theorem for a cover with disconnected intersection~\cite[\S1.2]{Hatcher} gives the following relations:
		\begin{itemize}
			\item from $S_{\rm out}$: $x_j'=x_j$ for $j>r$, and $x_1'\cdots x_r'=x_1\cdots x_r$;
			\item from $S_{\rm in}$: $\alpha_+\gamma\alpha_+^{-1}=c\,(\alpha_-\gamma\alpha_-^{-1})\,c^{-1}$ for each loop $\gamma$ in $S_{\rm in}$. Taking $\gamma$ to be the loop around puncture $i\le r$ gives $x_i'=c\,x_ic^{-1}$.
		\end{itemize}
		Since $x_i'=\psi_c(x_i)$ for all $i$, the first relation becomes $cwc^{-1}=w$, and
		\[
		\pi_1(Y)=\langle x_1,\dots,x_p,c\mid cwc^{-1}=w\rangle .
		\]
		The bottom face includes by $x_i\mapsto x_i$, and the top face includes by $x_i\mapsto\psi_c(x_i)$.
		
		The space $X\setminus N(C)$ is obtained from $Y$ by gluing the bottom face to the top face via $h$. This is a graph of spaces with one vertex space $Y$ and one edge space $D_p$, so $\pi_1(X\setminus N(C))$ is the HNN extension of $\pi_1(Y)$ with stable letter $\tau=[\nu]$ and the relations
		\[
		\tau\,\gamma\,\tau^{-1}=\psi_c(h_*\gamma)\qquad\bigl(\gamma\in\pi_1(D_p,*)\bigr),
		\]
		where $\gamma$ is included through the bottom face and $h_*\gamma$ through the top face~\cite[\S1.B]{Hatcher}. This presentation holds whether or not the face inclusions are injective on $\pi_1$.  Taking $\gamma=x_i$ gives the stated presentation. Every generator is a meridian by construction.
		\end{proof}

	\section{The augmented link $T(p,q)\cup A\cup C$}\label{sec:aug}

	Throughout this section $p\ge2$, $q\ge1$, $1\le r\le p-1$, and $g=\gcd(p,q)$, $p'=p/g$, $q'=q/g$; we allow $g=1$. We use the notation $\bu$, $U$, $\lambda_c$, $G_i$, $c(x)$, $N_i(k)$, and $S_c$ from the introduction.
	
	\subsection{Components and colors}\label{sec:colors}
	The permutation of $\delta$ is a $p$-cycle, so the permutation $j\mapsto j+q\pmod p$ of $\delta^q$ has the $g$ residue classes mod $g$ as its orbits, each of length $p'$. Hence $L=T(p,q)=\widehat{\delta^q}$ has $g$ components $L_0,\dots,L_{g-1}$, each a $T(p',q')$ torus knot, and strand $j$ lies on $L_{c(j)}$ with $c(j)\equiv(j-1)\bmod g$. We call $c(j)$ the \emph{color} of strand $j$, and write $\mathcal C_c=\{j: c(j)=c\}$. The full twist $\Delta_r^2$ is a pure braid, so inserting it does not change the components or their colors.
	
	Let $A$ be the braid axis and $C$ the circle around strands $1,\dots,r$ (Figure~\ref{fig:link}). Then
	\[
	\lk(L_c,A)=p',\qquad \lk(C,A)=0,\qquad \lk(L_c,C)=\lambda_c,
	\]
	the last because $C$ bounds a disc meeting $L_c$ in the $\lambda_c$ strands $j\le r$ of color $c$. The link $L\cup A\cup C$ has $g+2$ components: $L_0,\dots,L_{g-1}$, $A$, and $C$.
	
	\subsection{The formula}
	Let $S$ be the $p\times p$ cyclic shift matrix, $S_{j,j+1}=1$ with indices mod $p$, so that $(S^q)_{j,l}=1$ exactly when $l\equiv j+q\pmod p$. Let $D=\operatorname{diag}(d_1,\dots,d_p)$ with $d_j=v$ for $j\le r$ and $d_j=1$ for $j>r$, and let $\mathbf{1}$ be the all-ones vector. Set
	\[
	M=U^{q'}S^qD,\qquad \mu_c=U^{p'q'}v^{\lambda_c}\quad(0\le c\le g-1).
	\]
	
	\begin{theorem}\label{thm:aug}
		Let $a$ be the variable of $A$, and let
		\[
		S_c(\bu,v,a)=\sum_{n=0}^{p'-1}\ \sum_{\substack{0\le i\le r-1\\ i\equiv c\ (\mathrm{mod}\ g)}} a^{p'-1-n}\,U^{nq'}\,G_i\,v^{N_i(n)} ,
		\]
		so that $S_c(\bu,v,1)=S_c(\bu,v)$. Then
		\[
		\Delta_{L\cup A\cup C}(\bu,v,a)\;\ueq\;\sum_{c=0}^{g-1}(1-u_c)\prod_{c'\ne c}\bigl(a^{p'}-\mu_{c'}\bigr)\,S_c(\bu,v,a).
		\]
	\end{theorem}
	
	\begin{proof}
		 By Lemma~\ref{lem:presentation} with $\beta=\delta^q$, the group of $L\cup A\cup C$ has generators $x_1,\dots,x_p,c,\tau$ and relations
		\[
		R_i=\tau x_i\tau^{-1}\,\psi_c\bigl(\delta^q(x_i)\bigr)^{-1}\quad(1\le i\le p),\qquad R_c=c\,w\,c^{-1}w^{-1},\quad w=x_1\cdots x_r .
		\]
		We abelianize by $x_j\mapsto u_{c(j)}$, $c\mapsto v$, $\tau\mapsto a$. This is the abelianization with the meridian basis of $H_1$. Indeed, $R_c$ is a commutator. By Equation \ref{eq:delta}, $\delta(x_i) = x_1 x_{i+1} x_1^{-1}$ for $i$ indexed mod $p$. For $i=p$ this gives $x_1x_1x_1^{-1}=x_1$. Let $W_1 = x_1$ and $W_q = \delta(W_{q-1})x_1$. Then for all $i$ taken mod $p$, \begin{equation}\label{eq:conj}\delta^q(x_i)=\delta(W_{q-1})x_1x_{i+q}x_1^{-1}\delta(W_{q-1})^{-1}= W_qx_{i+q}W_q^{-1}.\end{equation} Thus $\delta^q(x_i)$ is conjugate to $x_{i+q}$, and $\psi_c$ conjugates by definition, so in $H_1$ the relation $R_i$ says $[x_i]=[x_{i+q}]$ (indices mod $p$). Adding $q$ modulo $p$ connects exactly the indices in one residue class mod $g$, so $H_1\cong\mathbb Z^{g+2}$ is free on the classes of meridians of $L_0,\dots,L_{g-1}$, $C$, and $A$.
		
		Let $B=J(\delta^q)$ be the abelianized Jacobian of $\delta^q$ for the above abelianization (the \emph{colored Burau matrix}). Since $\psi_c$ multiplies the derivative with respect to $x_j$ by $c$ for $j\le r$, we have $J(\psi_c)=D$. As per Equation \ref{eq:conj}, $\delta^q$ and $\psi_c$ preserve the abelianization, so the chain rule of Section~\ref{sec:fox} gives $J(\psi_c\circ\delta^q)=BD$. Hence, by Lemma~\ref{lem:alexjac}(b) with $\phi=\psi_c\circ\delta^q$ and $y_1=c$,
		\[
		\ab\Fox{R_i}{x_j}=a\delta_{ij}-(BD)_{ij},\qquad
		\ab\Fox{R_i}{\tau}=1-u_{c(i)},
		\]
		\[
		\ab\Fox{R_c}{x_j}=\begin{cases}(v-1)\,G_{j-1},&j\le r,\\0,&j>r,\end{cases}\qquad
		\ab\Fox{R_c}{\tau}=0.
		\]
		For the third entry, $\Fox{R_c}{x_j}=c\,\Fox{w}{x_j}-cwc^{-1}w^{-1}\Fox{w}{x_j}$ and $\ab\Fox{w}{x_j}=\ab(x_1\cdots x_{j-1})=G_{j-1}$ for $j\le r$. Let $\boldsymbol\tau$ be the column vector with entries $\tau_i=1-u_{c(i)}$, and $\boldsymbol\gamma=(G_0,G_1,\dots,G_{r-1},0,\dots,0)^{\mathsf T}$.
		
		The presentation of the group of $L\cup A\cup C$ has deficiency one, every generator is a meridian, and $L\cup A\cup C$ has $g+2\ge3$ components. By Lemma~\ref{lem:FoxTorres}, deleting the $c$-column gives
		\begin{equation}\label{eq:2}
		(v-1)\,\Delta_{L\cup A\cup C}\ueq\det Q,\qquad
		Q=\begin{pmatrix} aI-BD & \boldsymbol\tau\\ (v-1)\,\gamma^{\mathsf T} & 0\end{pmatrix},
		\end{equation}
		with rows $R_1,\dots,R_p,R_c$ and columns $x_1,\dots,x_p,\tau$.
		
		Since $\delta(x_j)$ is conjugate to $x_{j+1}$, induction gives $\ab(W_q)=\prod_{m=1}^{q}u_{c(m)}$. As $m$ runs over $1,\dots,q$, the color $c(m)$ runs through the residues mod $g$ periodically, each exactly $q'$ times, so $\ab(W_q)=U^{q'}$. Again, $\delta^q$ preserves the colored abelianization; differentiating the conjugate $\delta^q(x_i)=W_qx_{i+q}W_q^{-1}$,
		\[
		\Fox{\delta^q(x_i)}{x_j}=\Fox{W_q}{x_j}+W_q\,[\,j\equiv i+q\,]-W_qx_{i+q}W_q^{-1}\Fox{W_q}{x_j},
		\]
		and abelianizing gives
		\begin{align*}
		(B)_{i,j}&=U^{q'}[j\equiv i+q]+(1-u_{c(i)})\ab \Fox{W_q}{x_j}\\
		B &= U^{q'}S^q+\boldsymbol{\tau}\omega^{\mathsf T}, \omega_j = \ab \Fox{W_q}{x_j}.
		\end{align*}
		That is, $B$ is a weighted permutation matrix plus a rank-one matrix whose column factor is $\boldsymbol\tau$. Hence $BD=M+\boldsymbol\tau\,(D\omega)^{\mathsf T}$.
		
		For each $j$, add $(D\omega)_j$ times the $\tau$-column of $Q$ to its $x_j$-column. In row $i\le p$ this adds $\tau_i(D\omega)_j$, which cancels the $(i,j)$ entry of $-\boldsymbol\tau(D\omega)^{\mathsf T}$. The last row is unchanged, since its $\tau$-entry is $0$. These column operations do not change the determinant, and by the identity $\det\begin{pmatrix}X&b\\c^{\mathsf T}&0\end{pmatrix}=-c^{\mathsf T}\adj(X)\,b$,
		\begin{equation}\label{eq:4}
		\det Q=\det\begin{pmatrix} aI-M & \boldsymbol\tau\\ (v-1)\boldsymbol\gamma^{\mathsf T} & 0\end{pmatrix}=-(v-1)\,\boldsymbol\gamma^{\mathsf T}\adj(aI-M)\,\boldsymbol\tau .
		\end{equation}
		
		Since $(M\boldsymbol{y})_j=U^{q'}d_{j+q}\,y_{j+q}$, the matrix $M$ is monomial with underlying permutation $j\mapsto j+q\pmod p$, whose cycles are the color classes $\mathcal C_c$, each of length $p'$. Since $M$ never mixes strands of different colors, listing the strands color class by color class makes it block diagonal: $M=\bigoplus_c M_c$, where $\bigoplus$ denotes the block-diagonal sum and $M_c$ is the $p'\times p'$ block of $M$ on $\mathcal C_c$. Going once around $\mathcal C_c$ multiplies the weights $U^{q'}d_j$ for $j\in\mathcal C_c$, so $M_c^{p'}=\mu_cI$ with $\mu_c=U^{p'q'}v^{\lambda_c}$. Hence
		\[
		(aI-M_c)\sum_{n=0}^{p'-1}a^{p'-1-n}M_c^n=(a^{p'}-\mu_c)I ,
		\]
		so $\det(aI-M_c)=a^{p'}-\mu_c$  and $\adj(aI-M_c)=\sum_{n=0}^{p'-1}a^{p'-1-n}M_c^n$. For a block-diagonal matrix,
		\[
		\adj\Bigl(\bigoplus_c X_c\Bigr)=\bigoplus_c\bigg[\Bigl(\prod_{c'\ne c}\det X_{c'}\Bigr)\adj X_c\bigg],
		\]
		so
		\[
		\adj(aI-M)=\bigoplus_{c}\bigg[\Bigl(\prod_{c'\ne c}(a^{p'}-\mu_{c'})\Bigr)\sum_{n=0}^{p'-1}a^{p'-1-n}M_c^{n}\bigg].
		\]
		
		On $\mathcal C_c$ the vector $\boldsymbol\tau$ is the constant $1-u_c$. By induction on $n$, $(M^n\mathbf{1})_j=U^{nq'}\prod_{m=1}^n d_{j+mq}$, and the indices $j+mq$ stay in $\mathcal C_{c(j)}$. Take $j=i+1\le r$. Then $\gamma_j=G_i$, $c(j)\equiv i\bmod g$, and $d_{i+1+mq}=v$ exactly when $(i+1+mq)\bmod p\in\{1,\dots,r\}$, that is, when $(mq+i)\bmod p<r$. So $\prod_{m=1}^n d_{i+1+mq}=v^{N_i(n)}$. Collecting terms by color,
		\[
		\gamma^{\mathsf T}\adj(aI-M)\,\boldsymbol\tau=\sum_{c}(1-u_c)\prod_{c'\ne c}(a^{p'}-\mu_{c'})\sum_{n=0}^{p'-1}\sum_{\substack{0\le i<r\\ i\equiv c\,(g)}}a^{p'-1-n}\,U^{nq'}G_i\,v^{N_i(n)}.
		\]
		Combining with Equations \ref{eq:2} and \ref{eq:4} proves the theorem.
	\end{proof}

	\subsection{Proof of Theorem~\ref{thm:A}}
	
	\begin{proof}[Proof of Theorem~\ref{thm:A}]
		Let $g=\gcd(p,q)$. The link $T(p,q)\cup A\cup C$ has $g+2\ge3$ components, and $\lk(L_c,A)=p'$ for all $c$, while $\lk(C,A)=0$. By the Torres condition~\eqref{eq:torres} applied to $\ell=A$,
		\[
		\Delta_{T(p,q)\cup A\cup C}(\bu,v,1)\ueq(U^{p'}-1)\,\Delta_{T(p,q)\cup C}(\bu,v).
		\]
		Setting $a=1$ in Theorem~\ref{thm:aug}, \begin{align*}
			\Delta_{T(p,q)\cup C}(\bu,v) & \ueq \dfrac{\displaystyle\sum_{c=0}^{g-1}(1-u_c)\displaystyle\prod_{c'\neq c}(1-\mu_{c'})S_c(\bu,v)}{U^{p'}-1}\\
			& \ueq \dfrac{\displaystyle\sum_{c=0}^{g-1}(1-u_c)\displaystyle\prod_{c'\neq c}(1-U^{p'q'}v^{\lambda_c})S_c(\bu,v)}{U^{p'}-1}.
		\end{align*} When $g=1$, the produce over $c'\neq c$ is empty, $U = u$, $G_i = u^i$ and $k$ runs up to $p-1$ giving $$S_c(u,v)\ueq\sum_{k=0}^{p-1}\sum_{i=0}^{r-1}u^{kq+i}v^{N_i(k)}.$$We then have by Theorem~\ref{thm:aug} $$\Delta_{T(p,q)\cup C}(u,v)\ueq \dfrac{1-u}{1-u^p}\sum_{k=0}^{p-1}\sum_{i=0}^{r-1}u^{kq+i}v^{N_i(k)}.$$
	\end{proof}

	\section{The Alexander polynomials of twisted torus links}\label{sec:links}

	In this section we prove the formulae for the Alexander polynomials, first for twisted torus links and then for knots. 
	
	\begin{proof}[Proof of Theorem~\ref{thm:links}]
		Let $g = \gcd(p,q)\geq 2$. By Theorem~\ref{thm:A},
		\[
		\Delta_{T(p,q)\cup C}(\bu,v)\ueq\frac{\sum_c(1-u_c)\prod_{c'\ne c}\bigl(1-U^{p'q'}v^{\lambda_{c'}}\bigr)\,S_c(\bu,v)}{U^{p'}-1}.
		\]
		The link $T(p,q)\cup C$ has $d=g+1\ge3$ components. We then apply~\eqref{eq:SW} in the case $d\ge3$, with $L_d=C$, $\lambda_c=\lk(L_c,C)$, and $n=-s$. A $-1/s$ surgery on $C$ inserts $(\theta_r)^s$, so the resulting link is $T(p,q;r,s)$ with unchanged colors. Thus
		\[
		\Delta_{L\cup C}(\bu,V^{s})\ueq(V-1)\,\Delta_{T(p,q;r,s)}(\bu).
		\]
		Since $r\ge1$, some $\lambda_c\ge1$, so $V=V=u_0^{\lambda_0}\cdots u_{g-1}^{\lambda_{g-1}}\ne1$ and we may divide by $V-1$.
	\end{proof}
	
	We now prove the equivalent result for twisted torus knots, giving an explicit formula for $\Delta_{T(p,q;r,s)}$ for all $p,q,r,s$ with $\gcd(p,q)=1$.
	
		\begin{proof}[Proof of Theorem~\ref{thm:knots}]
		Let $\gcd(p,q)=1$. A $1/-s$ on $C$ yields the knot $T(p,q;r,s)$. By the second part of Theorem \ref{thm:aug}, $$\Delta_{T(p,q)\cup C}(u,v)\ueq \dfrac{1-u}{1-u^p}\sum_{k=0}^{p-1}\sum_{i=0}^{r-1}u^{kq+i}v^{N_i(k)}.$$ Applying Equation \ref{eq:SW} with $\ell=C$, $\lambda = \lk(T(p,q),C)=r$
		\begin{align*}\Delta_{T(p,q;r,s)}(t)& \ueq\dfrac{1-t}{1-t^r}\Delta_{T(p,q)
		\cup C}(t,t^{rs})\\
	& \ueq \dfrac{(1-t)^2}{(1-t^p)(1-t^r)}\sum_{k=0}^{p-1}\sum_{i=0}^{r-1}t^{kq+i+rs\,N_i(k)}.\end{align*}
	\end{proof}

	\section{Degree, fiberedness, and non-L-space knots}\label{sec:degree}

		\subsection{The counting function}\label{sec:N}
	
	Throughout, $\gcd(p,q)=1$ and $1\le r\le p-1$. Let $q^*$ denote the inverse of $q$ modulo $p$.
	
	\begin{lemma}\label{lem:Nforms}
		For $0\le i\le r-1$ and $0\le k\le p-1$:
		\begin{enumerate}
			\item[(a)] $N_i(0)=0$ and $N_i(k)=N_i(k-1)+[(kq+i)\bmod p<r]$;
			\item[(b)] $N_i(k)=\sum_{m=1}^{k}\Bigl(\bigl\lfloor\tfrac{mq+i}{p}\bigr\rfloor-\bigl\lfloor\tfrac{mq+i-r}{p}\bigr\rfloor\Bigr)$;
			\item[(c)] $N_i(k)=\#\{\rho\in\{0,\dots,r-1\}:\ 1\le (\rho-i)q^*\bmod p\le k\}$.
		\end{enumerate}
		In particular, $N_i(k)$ depends only on $q\bmod p$, it is nondecreasing in $k$ with steps $0$ or $1$, and $N_i(p-1)=r-1$.
	\end{lemma}
	
	\begin{proof}
		Part (a) is the definition. For (b), write $j=a p+b$ with $0\le b<p$. Then $\lfloor j/p\rfloor-\lfloor (j-r)/p\rfloor$ is $1$ if $b<r$ and $0$ otherwise. For (c), note that $(mq+i)\bmod p=\rho$ if and only if $m\equiv(\rho-i)q^*\pmod p$, and each residue class contains at most one $m$ with $1\le m\le k<p$. Finally, as $m$ runs over $1,\dots,p-1$, the residue $(mq+i)\bmod p$ takes every value except $i$. Exactly $r-1$ of these values lie in $\{0,\dots,r-1\}$.
	\end{proof}
	
	\begin{lemma}[Pairing]\label{lem:pairing}
		For all $i$ and $k$, $\;N_i(k)+N_{r-1-i}(p-1-k)=r-1$.
	\end{lemma}
	
	\begin{proof}
		Substitute $m=p-m'$ in the definition of $N_{r-1-i}(p-1-k)$, so that $k+1\le m'\le p-1$. Let $y\equiv(m'q+i)\bmod p$. Then $(mq+r-1-i)\bmod p=(r-1-y)\bmod p$, and this is less than $r$ exactly when $y\le r-1$. Hence $N_{r-1-i}(p-1-k)=N_i(p-1)-N_i(k)=r-1-N_i(k)$ by Lemma~\ref{lem:Nforms}.
	\end{proof}

	\subsection{Non-fibered knots}
	
	\begin{proof}[Proof of Theorem~\ref{thm:degree}]
		Recall $E(k,i)=kq+i+rsN_i(k)$. Let $F_{p,q,r,s}(t)=\displaystyle\sum_{k=0}^{p-1}\sum_{i=0}^{r-1}t^{E(k,i)}.$ By Theorem~\ref{thm:knots}, $$(1-t^p)(1-t^r)\Delta_{T(p,q;r,s)}(t)\ueq(1-t)^2F_{p,q,r,s}(t).$$Span is additive under multiplication in $\mathbb Z[t^{\pm1}]$, so $\spn\Delta_{T(p,q;r,s)}=\spn F_{p,q,r,s}-(p+r-2)$. Since $F_{p,q,r,s}$ is a sum of monomials with coefficient $+1$, no cancellation occurs, and $\spn F_{p,q,r,s}=\max E-\min E$. By Lemma~\ref{lem:pairing},
		\begin{equation}\label{eq:sides}
		E(k,i)+E(p-1-k,r-1-i)=(p-1)q+(r-1)+rs(r-1)=:c_{p,q,r,s}.
		\end{equation}
		Hence $\min E=c_{p,q,r,s}-\max E$, which completes the proof.\end{proof}
	
	\begin{proof}[Proof of Corollary~\ref{cor:fiber}]
		Compare the top coefficients in $(1-t^p)(1-t^r)\Delta_{T(p,q;r,s)}\ueq(1-t)^2F_{p,q,r,s}$ of Equation \ref{eq:sides}. On the left-hand side, the top coefficient is $\pm$ the leading coefficient of $\Delta$. On the right-hand side of \ref{eq:sides}, it is the top coefficient of $F$, namely $\#\{(k,i):E(k,i)=\max E\}$. The Alexander polynomial of a fibered knot is monic~\cite{Seifert}. 
	\end{proof}
	
	\subsection{Non-L-Space Knots} By Ozsv\'ath and Szab\'o \cite{OzsvathSzabo}, the Alexander polynomial of an L-space knot must have all coefficients $\pm 1$. Due to the nature of the formula in Theorem \ref{thm:knots}, it is then not hard to produce knots which are not L-space knots. A non-monic Alexander polynomial as in \ref{cor:fiber} already rules out an L-space knot. The following lemma detects non-L-space knots whose polynomial is monic. 
	
	\begin{lemma}\label{lem:m0m1}
		Let $gcd(p,q)=1, 2\leq r \leq p-1$ and $s\in \mathbb Z$ be arbitrary. Let $m_0$ be the number of pairs $(k,i)$ at which $E(k,i)$ is maximal, and $m_1$ the number at which $E(k,i) = \max E -1$. Then, up to a unit,$$\Delta_{T(p,q;r,s)}(t) = m_0t^D + (m_1-2m_0)t^{D-1} + \cdots.$$In particular, if $m_0=1$ and $m_1=0$, then $\Delta_{T(p,q;r,s)}$ is monic and $T(p,q;r,s)$ is not an L-space knot.
	\end{lemma}
	
	\begin{proof}
		As in the proof of Theorem \ref{thm:degree}, let $F_{p,q,r,s}(t)=\displaystyle\sum_{k=0}^{p-1}\sum_{i=0}^{r-1}t^{E(k,i)}.$ Setting $x = t^{-1}$ we have $$t^{-\max E}F_{p,q,r,s}(t)=\sum_{k,i}x^{\max E - E(k,i)} = m_0 + m_1x + \cdots.$$ Letting $G(x) = m_0 + m_1x + \cdots$ as above and $H(x) = \dfrac{(1-x)^2}{((1-x^p)(1-x^r))}$, $\Delta_{T(p,q;r,s)}(x) \ueq G(x)H(x)$. Written as a product of power series, $$H(x) = (1-2x+x^2)(1+x^p + \cdots)(1+x^r + \cdots).$$ As $p,r\geq 2$, we have $$\Delta_{T(p,q;r,s)}(x) \ueq (m_0 + m_1x + \cdots)(1-2x + \cdots) = m_0 + (m_1-2m_0)x+\cdots.$$ It follows then that $$\Delta_{T(p,q;r,s)}(t) = m_0t^D + (m_1-2m_0)t^{D-1} + \cdots, \quad D = \text{deg }\Delta_{T(p,q;r,s)}.$$ If $m_0 =1 $ and $m_1=0$, the second coefficient is $-2$; by \cite{OzsvathSzabo}, $T(p,q;r,s)$ is not an L-space knot.
	\end{proof}
	
	\begin{proof}[Proof of Theorem~\ref{thm:lspace}] As $r<q<p$, we have $T(p,q;r,s) = T(q,p;r,s)$ by Lee \cite[Lemma 5]{Lee}. We apply Theorem \ref{thm:knots} to $T(q,p;r,s)$ (note the proof does not rely on $q<p$). We write $r|s| = p+1+d$ with $d\geq 2$. Since $p \equiv -1 \mod q$, we have $(mp+i)\mod q \equiv (i-m)\mod q$. For $1\leq m \leq i$ this equals $i-m <r$, so $m$ is counted in $N_i(k)$. For $i < m \leq q-1$ it equals $q+i-m$, which is less than $r$ exactly when $m > q+i-r$. Hence $$N_i(k) = \min(k,i) + \max(0,k-(q+i-r)),$$ and $E(k,i) = kp+i - (p+1+d)N_i(k).$ Then $k$ falls in three regions:\begin{enumerate}
			\item[(I)] if $k\leq i$, then $E(k,i)=i-k(1+d) \leq r-1$;
			\item[(II)] if $i\leq k \leq q+i-r$, then $E(k,i)=kp-i(p+d)\leq (q-r)p-id$, with equality exactly when $k=q+i-r$;
			\item[(III)] if $k>q+i-r$, then $E(k,i)=(p+1+d)(q-r) +i-k(1+d) \leq (q-r)p-(1+d)-id$. 
		\end{enumerate} Thus $\max E = (q-r)p$ attained only at $(k,i) = (q-r,0)$. Hence $m_0=1$.
		
		In region (II), the values of $E(k,i)$ with $i\geq 1$ are at most $\max E -d$, and those with $i=0$ and $k<q-r$ are at most $\max E - p$. In region (III), $E(k,i)$ is at most $\max E -3$, and in region (I) at most $r-1 <(q-r)p-1$. So $m_1 = \max E-1 = 0$; by Lemma \ref{lem:m0m1}, $\Delta_{T(p,q;r,s)}(t) = t^D -2t^{D-1} + \cdots$. As desired we have $T(p,q;r,s)$ has a monic Alexander polynomial yet is not an L-space knot by \cite{OzsvathSzabo}.\end{proof}
	
	
	The below edge case was first identified by the author's student Ethan Elliott with the aid of chatGPT. The proof of correctness of the coefficients follows closely that of Theorem \ref{thm:lspace}.
	
	\begin{proposition}\label{prop:ethan}
		For $p\ge5$, let $K_p=T(p,\,p-1;\,3,\,-p)$. Up to a unit, $\Delta_{K_p}$ can be written as a polynomial with nonzero constant term:
		\[
		\Delta_{K_p}(t)=t^{D_p}-2\,t^{D_p-1}+\cdots-2t+1,\qquad D_p=p^2-p+2 .
		\]
		In particular, $\Delta_{K_p}$ is monic, but $K_p$ is not an L-space knot, and the knots $K_p$ are pairwise distinct.
	\end{proposition}
	
	\begin{proof}
		Here $q=p-1$, $r=3$, and $rs=-3p$. Since $q\equiv-1\pmod p$, we have $(mq+i)\bmod p=(i-m)\bmod p$. For $1\le m\le i$ this is $i-m<3$; for $i<m\le p-1$ it is $p+i-m$, which is less than $3$ exactly when $m\ge p-2+i$. Hence, for $0\le i\le2$ and $0\le k\le p-1$,
		\[
		N_i(k)=\min(k,i)+\max\bigl(0,\;k-(p-3+i)\bigr),
		\]
		and $E(k,i)=k(p-1)+i-3p\,N_i(k)$. Put $M=(p-3)(p-1)$. There are three ranges:
		\begin{enumerate}
			\item[(I)] if $k\le i$, then $N_i(k)=k$ and $E(k,i)=i-(2p+1)k\le2$;
			\item[(II)] if $i\le k\le p-3+i$, then $N_i(k)=i$ and $E(k,i)=(p-1)k-(3p-1)i\le M-2p\,i$, with equality exactly at $k=p-3+i$;
			\item[(III)] if $k\ge p-3+i$, then $N_i(k)=k-p+3$ and $E(k,i)=3p(p-3)+i-(2p+1)k\le M-2p\,i$, with equality exactly at $k=p-3+i$.
		\end{enumerate}
		So $\max E=M$, attained only at $(k,i)=(p-3,0)$. Every other value is at most $M-2$. Indeed:
		\begin{itemize}
			\item for $i=0$ and $k\ne p-3$, ranges (II) and (III) give at most $M-(p-1)$;
			\item for $i\ge1$ they give at most $M-2p$;
			\item range (I) gives at most $2<M-1$, since $M\ge8$.
		\end{itemize}
		Hence $F_{p,q,r,s}(t)=t^{M}+(\text{terms of degree}\le M-2)$.
		
		By Theorem~\ref{thm:knots},
		\begin{equation}\label{eq:str}
		(1-t^p)(1-t^3)\,\Delta_{K_p}(t)\ueq(1-2t+t^2)\,F_{p,q,r,s}(t).
		\end{equation}
		The two highest terms of the right side are then $t^{M+2}-2t^{M+1}$. The left hand side of Equation \ref{eq:str} has $(1-t^p)(1-t^3)=t^{p+3}-t^p-t^3+1$, and its next term after $t^{p+3}$ lies three degrees lower. So the two highest coefficients of the left side are the two highest coefficients of $\Delta_{K_p}$. Hence these are $1$ and $-2$, up to a common sign.
		
		For the degree, Theorem~\ref{thm:degree} with $\max E=M$ and $(p-1)q+(r-1)(1+rs)=(p-1)^2+2(1-3p)$ gives
		\[
		\spn\Delta_{K_p}=2(p-3)(p-1)-(p-1)^2-2(1-3p)-(p+1)=p^2-p+2 .
		\]
		By Ozsv\'ath and Szab\'o~\cite{OzsvathSzabo}, the nonzero coefficients of the Alexander polynomial of an L-space knot are $\pm1$, so the coefficient $-2$ shows that $K_p$ is not an L-space knot. The spans are distinct, so the knots are pairwise distinct.
	\end{proof}
	
	\begin{remark}
		The hypothesis $p\ge5$ is needed. For $p=4$ the term at $(k,i)=(0,2)$ has $E=2=M-1$, and $\Delta_{T(4,3;3,-4)}=\Delta_{T(3,8)}$ satisfies the Ozsv\'ath--Szab\'o condition. The same argument applies to $T(p,p-1;r,s)$ whenever $2\le r\le p-2$ and $r|s|\ge p+2$.
	\end{remark}

	\section{Questions}\label{sec:questions}

	\begin{question}
		Determine exactly when $\max E$ is attained more than once, in terms of $p$, $q\bmod p$, $r$, and $s$. This would give a closed-form non-fiberedness criterion for negatively twisted torus knots.
	\end{question}

	\begin{question}
		For twisted torus links, the numerator in Theorem~\ref{thm:links} has coefficients of both signs, so the argument of Theorem~\ref{thm:degree} does not apply directly. Determine the Newton polytope of $\Delta_{T(p,q;r,s)}$, equivalently its Alexander norm ball~\cite{McMullen}, and give a fiberedness criterion analogous to Corollary~\ref{cor:fiber}.
	\end{question}

	\subsection*{Acknowledgments}
	The author was supported by the Elizabeth Gaines Mann Foundation.
	
	\subsection*{Tool Use.} The main results, including the closed formula, are the author's, and all proofs were written by the author. The family of twisted torus knots in Proposition \ref{prop:ethan} were found and verified by Ethan Elliott with the assistance of ChatGPT (GPT-6.1-Sol) and Claude (Opus-5.5).


\begin{thebibliography}{99}

		
		\bibitem{Dean} J. Dean, \emph{Small Seifert-fibered Dehn surgery on hyperbolic knots}, Algebr. Geom. Topol. \textbf{3}, no.1 (2003), 435--472.
		
		\bibitem{Morton} H. R. Morton, \emph{The Alexander polynomial of a torus knot with twists}, J. Knot Theory Ramifications \textbf{15} (2006), 1037--1047.
		
		\bibitem{AdnanPark} Adnan and K. Park, \emph{The Alexander polynomial of twisted torus knots}, J. Knot Theory Ramifications \textbf{34}, no. 08 (2025), p.2550038.
		
		\bibitem{SilverWilliams} D. S. Silver and S. G. Williams, \emph{Mahler measure of Alexander polynomials}, J. London Math. Soc. (2) \textbf{69} (2004), 767--782.
		
		\bibitem{AdnanPark2} Adnan and K. Park, \emph{Infinite families of non-fibered twisted torus knots}, Korean Journal of Mathematics, \textbf{34}, no. 2 (2026), 233-242.
		
		\bibitem{OzsvathSzabo} P. Ozsv\'ath and Z. Szab\'o, \emph{On knot Floer homology and lens space surgeries}, Topology \textbf{44} (2005), 1281--1300.
		
		\bibitem{Vafaee} F. Vafaee, \emph{On the knot Floer homology of twisted torus knots}, Int. Math. Res. Not. IMRN \textbf{2015}, no.~15, 6516--6537.
		
			\bibitem{Fox53} R. H. Fox, \emph{Free differential calculus. I. Derivation in the free group ring}, Ann. of Math.  \textbf{57}, no. 3 (1953), 547--560.
		
		\bibitem{CrowellFox} R. H. Crowell and R. H. Fox, \emph{Introduction to Knot Theory}, Grad. Texts in Math. \textbf{57}, Springer, 1977.
		
		\bibitem{Torres} G. Torres, \emph{On the Alexander polynomial}, Ann. of Math. (2) \textbf{57} (1953), 57--89.
		
		\bibitem{Kawauchi} A. Kawauchi, \emph{A Survey of Knot Theory}, Birkh\"auser, Basel, 1996.
		
		\bibitem{KasselTuraev} C. Kassel and V. Turaev, \emph{Braid Groups}, Grad. Texts in Math. \textbf{247}, Springer, 2008.
		
		\bibitem{Birman} J. S. Birman, \emph{Braids, Links, and Mapping Class Groups}, Ann. of Math. Stud. \textbf{82}, Princeton Univ. Press, 1974.
		
		
	
		
		\bibitem{Hatcher} A. Hatcher, \emph{Algebraic Topology}, Cambridge Univ. Press, 2002.
		
		
		
		\bibitem{Seifert} D. Rolfsen, \emph{Knots and Links}, No. 346 American Mathematical Soc. 2003.
		
		\bibitem{Lee} S. Lee, \emph{Knot types of twisted torus knots}, J. Knot Theory Ramifications, \textbf{26}, no. 12 (2017).
		
		\bibitem{McMullen} C. T. McMullen, \emph{The Alexander polynomial of a 3-manifold and the Thurston norm on cohomology}, Ann. Sci. \'Ec. Norm. Sup\'er. (4) \textbf{35} (2002), 153--171.
		

		

		
		
	\end{thebibliography}
\end{document}